\documentclass[amstex,12pt,russian,a4paper]{article}

\usepackage{mathtext}
\usepackage[cp1251]{inputenc}
\usepackage[T2A]{fontenc}
\usepackage[russian]{babel}
\usepackage[dvips]{graphicx}
\usepackage{amsmath}
\usepackage{amssymb}
\usepackage{amsxtra}
\usepackage{latexsym}
\usepackage{ifthen}

\renewcommand{\baselinestretch}{1.3}

\begin{document}

\def\Xint#1{\mathchoice
   {\XXint\displaystyle\textstyle{#1}}%
   {\XXint\textstyle\scriptstyle{#1}}%
   {\XXint\scriptstyle\scriptscriptstyle{#1}}%
   {\XXint\scriptscriptstyle\scriptscriptstyle{#1}}%
   \!\int}
\def\XXint#1#2#3{{\setbox0=\hbox{$#1{#2#3}{\int}$}
     \vcenter{\hbox{$#2#3$}}\kern-.5\wd0}}
\def\dashint{\Xint-}

\renewcommand{\abstractname}{}

\noindent
\renewcommand{\baselinestretch}{1.1}
\normalsize \large \noindent

\title{Геодезические касательные пространства к метрическим пространствам}

\author{В. В. Билет}

\medskip

\small {УДК 517.5}

\medskip

{\bf В. В. Билет} \small {(\textsf{Институт прикладной математики и механики
НАН Украины, г. Славянск})}

{\bf Р. Р. Салимов} \small {(\textsf{\textsf{Институт математики
НАН Украины, г. Киев}})}

\medskip

{\bf В. В. Білет} \small {(\textsf{Інститут прикладної математики і механіки НАН України, м. Слов'янськ})}

{\bf Р. Р. Салімов} \small {(\textsf{Інститут математики
НАН України, м. Київ})}

\medskip

{\bf V. V. Bilet} \small {(\textsf{Institute of Applied Mathematics and Mechanics of NASU,
Sloviansk})}

{\bf R. R. Salimov} \small {(\textsf{Institute of Mathematics of NASU,
Kyiv})}

\begin{abstract}

\textbf{Оценка площади образа круга для классов Соболева}

\medskip

Для регулярных гомеоморфизмов класса Соболева $W^{1,1}_{\textrm{loc}}$, обладающих $N$-свойством Лузина, установлена оценка площади образа круга в терминах угловой дилатации. Как следствие, получен аналог известной леммы Икома-Шварца для таких отображений.

\medskip

\textbf{Ключевые слова:} \emph{угловая дилатация, изопериметрическое неравенство, класс Соболева}

\medskip

\textbf{Оцінка площі образу кола для класів Соболєва}

\medskip

Для регулярних гомеоморфізмів класу Соболєва $W^{1,1}_{\textrm{loc}}$, які задовольняють $N$-властивості Лузіна, встановлено оцінку площі образу кола в термінах кутової дилатації. Як наслідок, отримано аналог відомої леми Ікома-Шварца для таких відображень.

\medskip

\textbf{Ключові слова:} \emph{кутова дилатація, ізопериметрична нерівність, клас Соболєва}

\medskip

\textbf{The estimation of the area of a disk image for Sobolev classes}

\medskip

For regular homeomorphisms of Sobolev class $W^{1,1}_{\textrm{loc}}$ having the Luzin
$N$-property, it is established the estimation of the area of a disk image in terms of an angular
dilatation. As a corollary, the analog of the well-known Ikoma-Schwartz lemma for such mappings is
obtained.

\medskip


\end{abstract}

\newpage

{\bf 1. Введение.} В данной статье получены точные оценки искажения площади образа круга при регулярных гомеоморфизмах класса Соболева $W^{1, 1}_{\rm loc},$
обладающих $N$-свойством (Лузина).

 Отметим, что впервые оценка площади образа круга при
 квазиконформных отображениях встречается в монографии Лаврентьева М. А., см. \cite{L}.
 В монографии \cite{BGMR}, см. предложение~3.7, получено уточнение неравенства Лаврентьева
 в терминах угловой дилатации. Также ранее в работах  \cite{LS} и \cite{Salimov1} верхние
  оценки искажения площади образа круга были получены методом  модулей.

Пусть $G$ -- область в комплексной плоскости $\mathbb C,$ то есть связное, открытое подмножество $\mathbb C.$ Напомним, что отображение $f: G\to\mathbb C$
называется \emph{регулярным} \emph{в точке} $z_0\in G,$ если в этой точке $f$ имеет полный
дифференциал и его Якобиан $J_f=|f_z|^{2}-|f_{\bar{z}}|^{2}\ne 0$ (см., например, I.~1.6 в
\cite{LV}). Гомеоморфизм $f$ класса Соболева $W_{\rm{loc}}^{1, 1}$ называется
\emph{регулярным}, если $J_{f}>0$ почти всюду (п.в.).

Говорят, что гомеоморфизм $f: G\to\mathbb C$ обладает $N$-\emph{свойством} (Лузина), если для
любого множества $E\subset G$ из условия $|E|=0$ следует, что $|f(E)|=0.$

\medskip

{\bf 2. Вспомогательные результаты.}

Обозначим через $L(r)$ длину кривой $f(re^{i\theta}),$ $0\leqslant\theta\leqslant 2\pi,$ и через
$S(r)=|f(B_r)|$~-- площадь $f(B_r).$ Везде далее полагаем $$B_r=\{z\in\mathbb C: |z|\leqslant r\}, \quad \mathbb B=\{z\in\mathbb C: |z|\leqslant 1\}.$$

\medskip

{\bf Лемма 1.}\label{lem1}
{\it Пусть $f:  \mathbb B\to \mathbb B$ ---  регулярный гомеоморфизм класса Соболева  $W^{1, 1}_{\rm loc},$ обладающий $N$-свойством. Тогда при $p\geqslant 2$
\begin{equation}\label{eqlem1}
L^{p}(r)\leqslant\delta_{p}(r)S'(r)
\end{equation}
для п.в. $r\in [0, 1],$ где
\begin{equation}\label{delta}
\delta_{p}(r)=\left(\int\limits_{\gamma_r} D^{\frac{1}{p-1}}_{p}(z)\, |dz |\right)^{p-1},\quad \gamma_{r}=\{z\in\mathbb C: |z|=r\}.
\end{equation}}

\medskip

{\it Доказательство.}
Для п.в. $r\in [0, 1]$, имеем
\begin{equation}
L(r)=\int\limits_{0}^{2\pi}|f_{\theta}(re^{i\theta})|\,d\theta=\int\limits_{0}^{2\pi}\, D_{p}^{\frac{1}{p}}(re^{i\theta})\, J_{f}^{\frac{1}{p}}(re^{i\theta})\, r\,d\theta,
\end{equation}
где
\begin{equation}\label{eq2}
D_{p}(re^{i\theta})=\frac{|f_{\theta}(re^{i\theta})|^{p}}{r^pJ_f(re^{i\theta})}\  -
\end{equation}
$p$-угловая дилатация и $J_f$ обозначает Якобиан  отображения $f$.

Применяя неравенство Гельдера, имеем
\begin{equation}
L(r)\leqslant \left(\int\limits_{0}^{2\pi} J_{f}(re^{i\theta})\, r\, d\theta\right)^{\frac{1}{p}}\left(\int\limits_{0}^{2\pi}\, D^{\frac{1}{p-1}}_{p}(re^{i\theta})\, r\, d\theta\right)^{\frac{p-1}{p}}
\end{equation}
для п.в. $r\in[0, 1].$
Используя тот факт, что \begin{equation}S'(r)\,=\,\int\limits_{0}^{2\pi}\, J_{f}(re^{i\theta})\, r\,  d\theta,\end{equation} получаем
\begin{equation}
L^p(r)\leqslant \delta_{p}(r)S'(r)
\end{equation}
для п.в. $r\in[0, 1].$
Лемма~1 доказана.

\medskip

В следующей лемме получено дифференциальное неравенство для функции $S(r).$

\medskip

{\bf Лемма 2.}\label{lem2}
{\it Пусть $f:  \mathbb B\to \mathbb B$ ---  регулярный гомеоморфизм класса Соболева  $W^{1, 1}_{\rm loc},$ обладающий $N$-свойством. Тогда при $p\geqslant 2$
\begin{equation}\label{eqlem2}
S'(r)\geqslant (4\pi)^{\frac{p}{2}}\, \delta_{p}^{-1}(r)\, S^{\frac{p}{2}}(r)
\end{equation}
для п.в. $r\in [0, 1],$
где $\delta_{p}(r)$ определено равенством \eqref{delta}}.

\medskip

{\it Доказательство.}
Из неравенства \eqref{eqlem1} следует, что для п.в. $r\in [0, 1]$
\begin{equation}
S'(r)\geqslant \delta_{p}^{-1}(r)\, L^{p}(r).
\end{equation}
Используя, далее, известное изопериметрическое неравенство \begin{equation}L^{2}(r)\geqslant 4\pi S(r),\end{equation} получаем дифференциальное неравенство \eqref{eqlem2}. Лемма~2 доказана.

\medskip

{\bf Лемма 3.}\label{lem3}
{\it Пусть $f:  \mathbb B \to \mathbb B$ ---  регулярный гомеоморфизм класса Соболева  $W^{1, 1}_{\rm loc},$ обладающий $N$-свойством. Тогда при $p=2$ имеет место оценка
\begin{equation}\label{eq1lem3}
S(r_1)\leqslant S(r_2)\exp\left\{-4\pi\int\limits_{r_1}^{r_2}\frac{dr}{\delta_{2}(r)}\right\},
\end{equation}
а при $p>2$ ---
\begin{equation}\label{eq1lem3}
S^{\frac{2-p}{2}}(r_1)-S^{\frac{2-p}{2}}(r_2)\leqslant (4\pi)^{\frac{p}{2}}\frac{p-2}{2}\int\limits_{r_1}^{r_2}\frac{dr}{\delta_{p}(r)},
\end{equation}
где $0<r_1\leqslant r_2<1$ и $\delta_{p}(r)$ определено равенством \eqref{delta}}.

\medskip

{\it Доказательство.}
Воспользуемся дифференциальным неравенством \eqref{eqlem2}. Тогда для п.в. $r\in [0, 1]$
\begin{equation}
\frac{S'(r)}{S^{\frac{p}{2}}(r)}\,dr\geqslant (4\pi)^{\frac{p}{2}}\frac{dr}{\delta_{p}(r)}.
\end{equation}
При $p=2,$ имеем
\begin{equation}
\frac{S'(r)}{S(r)}\, dr\geqslant 4\pi\, \frac{dr}{\delta_{2}(r)}.
\end{equation}
Интегрируя обе части последнего неравенства по $r\in [r_1, r_2]$, получаем
\begin{equation}
\int\limits_{r_1}^{r_2}\left(\ln S(r)\right)^{'}dr\geqslant 4\pi \int\limits_{r_1}^{r_2}\frac{dr}{\delta_{2}(r)}.
\end{equation}
Заметим, что функция $g(r)=\ln S(r)$ является неубывающей, тогда
\begin{equation}
\int\limits_{r_1}^{r_2}(\ln S(r))^{'}dr=\int\limits_{r_1}^{r_2}g'(r)dr\leqslant g(r_2)-g(r_1)=\ln\frac{S(r_2)}{S(r_1)}
\end{equation}
(см., напр., теорему IV. 7.4  в \cite{Sa})
и
\begin{equation}
\ln\frac{S(r_2)}{S(r_1)}\geqslant 4\pi \int\limits_{r_1}^{r_2}\frac{dr}{\delta_{2}(r)},
\end{equation}
откуда
\begin{equation}
S(r_1)\leqslant S(r_2)\exp\left\{-4\pi\int\limits_{r_1}^{r_2}\frac{dr}{\delta_{2}(r)}\right\}.
\end{equation}

При $p>2,$ имеем
\begin{equation}
\int\limits_{r_1}^{r_2}\left(\frac{S^{\frac{2-p}{2}}(r)}{\frac{2-p}{2}}\right)^{'}dr\geqslant (4\pi)^{\frac{p}{2}}\int\limits_{r_1}^{r_2}\frac{dr}{\delta_{p}(r)}.
\end{equation}
Заметим, что функция $g_{p}(r)=\frac{S^{\frac{2-p}{2}}(r)}{\frac{2-p}{2}}$ является неубывающей, тогда
\begin{equation}
\int\limits_{r_1}^{r_2}\left(\frac{S^{\frac{2-p}{2}}(r)}{\frac{2-p}{2}}\right)^{'}dr=\int\limits_{r_1}^{r_2}g_{p}'(r)\leqslant g_{p}(r_2)-g_{p}(r_1)=\frac{2}{2-p}\left(S^{\frac{2-p}{2}}(r_2)-S^{\frac{2-p}{2}}(r_1)\right)
\end{equation}
(см., теорему IV. 7.4  в   \cite{Sa}).

Таким образом, получаем
\begin{equation}
S^{\frac{2-p}{2}}(r_1)-S^{\frac{2-p}{2}}(r_2)\leqslant (4\pi)^{\frac{p}{2}}\frac{p-2}{2}\int\limits_{r_1}^{r_2}\frac{dr}{\delta_{p}(r)}.
\end{equation}
Лемма~3 доказана.

\medskip

{\bf 3. Основные результаты.}

Ниже приведена теорема об  оценке площади образа круга.

\medskip

{\bf Теорема 1.}\label{Th1} {\it Пусть $f:  \mathbb B\to \mathbb B$ ---  регулярный гомеоморфизм класса Соболева  $W^{1, 1}_{\rm loc},$ обладающий $N$-свойством. Тогда при $p=2$
имеет место оценка
\begin{equation}\label{eqks*}
|f(B_r)|\leqslant  \pi \exp\left\{-4\pi\int\limits_{r}^{1}\frac{dt}{\delta_2(t)}\right\},
\end{equation}
а при $p>2$ ---

\begin{equation}\label{eqks**}
|f(B_r)|\leqslant  \pi\left(1+(2\pi)^{p-1}(p-2)\int\limits_{r}^{1}\frac{dt}{\delta_{p}(t)} \right)^{-\frac{2}{p-2}}\,,
\end{equation}
где

\begin{equation*}
\delta_{p}(r)=\left(\int\limits_{\gamma_r} D^{\frac{1}{p-1}}_{p}(z)\, |dz |\right)^{p-1},\quad \gamma_{r}=\{z\in\mathbb C: |z|=r\}.
\end{equation*}}

\medskip

{\it Доказательство.} Теорема~1 следует из леммы~3 при условии, что $S(1)~\leqslant~|f(\mathbb B)|=\pi$.

\medskip


Из теоремы~1 вытекает аналог известной леммы Икома-Шварца (см. теорему~2 в \cite{Ikoma}).

\medskip

{\bf Теорема 2.}\label{Th2}
{\it Пусть $f:  \mathbb B\to \mathbb B$ ---  регулярный гомеоморфизм класса Соболева  $W^{1, 1}_{\rm loc},$ обладающий $N$-свойством, $f(0)=0$. Тогда при $p=2$ имеет место оценка
\begin{equation}
\liminf_{z\to 0}|f(z)|\exp\left\{2\pi\int\limits_{|z|}^{1}\frac{dt}{\delta_2(t)}\right\}\leqslant 1,
\end{equation}
а при $p>2$ ---
\begin{equation}
\liminf_{z\to 0}|f(z)|\left(1+(2\pi)^{p-1}(p-2)\int\limits_{|z|}^{1}\frac{dt}{\delta_{p}(t)} \right)^{\frac{1}{p-2}}\leqslant 1.
\end{equation}}

\medskip

{\it Доказательство.}
Положим $l_f(r)=\min\limits_{|z|=r}|f(z)|$.
Тогда, учитывая, что $f(0)=0$, получаем
$\pi\,l_{f}^{2}(r)\leqslant |f(B_r)|$ и, следовательно,
\begin{equation}\label{eqks1.20}l_{f}(r)\leqslant \sqrt{\frac{|f(B_r)|}{\pi}}.\end{equation}
Таким образом, учитывая  неравенства  \eqref{eqks*} и  \eqref{eqks**}, имеем
$$\liminf\limits_{z\to 0}\frac{|f(z)|}{\mathcal{R}_p(|z|)}=\liminf\limits_{r\to0}\frac{l_{f}(r)}{\mathcal{R}_p(r)}
\leqslant\liminf\limits_{r\to0}\sqrt{\frac{|f(B_r)|}{\pi}}\cdot\frac{1}{\mathcal{R}_p(r)}\leqslant1,$$ где
$$\mathcal{R}_p(r)=        \left(1+(2\pi)^{p-1}(p-2)\int\limits_{r}^{1}\frac{dt}{\delta_{p}(t)} \right)^{-\frac{1}{p-2}} $$
при $p>2$  и
$$\mathcal{R}_p(r)=\exp\left\{-2\pi\int\limits_{r}^{1}\frac{dt}{\delta_2(t)}\right\}$$
при $p=2$. Теорема~2 доказана.

\medskip

{\bf Следствие 1.}
Пусть $f:  \mathbb B\to \mathbb B$ ---  регулярный гомеоморфизм класса Соболева  $W^{1, 1}_{\rm loc},$ обладающий $N$-свойством, $f(0)=0$. Тогда при $p>2$ имеет место оценка
\begin{equation*}
\liminf_{z\to 0}|f(z)|\left(\int\limits_{|z|}^{1}\frac{dt}{\delta_{p}(t)} \right)^{\frac{1}{p-2}}\leqslant (2\pi)^{1-p}(p-2)^{\frac{1}{2-p}}.
\end{equation*}


\bigskip

\textbf{Виктория Викторовна Билет}

Институт прикладной математики и механики НАН Украины

ул. Генерала Батюка, 19,

г. Славянск, Украина, 84116,

e-mail: biletvictoriya@mail.ru;

\medskip

\textbf{Руслан Радикович Салимов}

Институт математики НАН Украины

ул. Терещенковская, 3,

г. Киев-4, Украина, 01004,

e-mail: ruslan623@yandex.ru

\end{document}